\documentclass[a4paper,12pt]{article}
\usepackage[top=2.5cm,bottom=2.5cm,left=2.5cm,right=2.5cm]{geometry}
\usepackage{amsfonts}
\usepackage{latexsym,amsmath,amssymb,amsfonts,epsfig,psfrag,url,graphics,ifpdf,multicol}
\usepackage{color,colordvi,amscd,amsthm} 
\usepackage{tikz}
\usepackage[numbers,sort&compress]{natbib}
\usepackage{indentfirst,graphics,epsfig}
\usepackage{graphicx}
\usepackage{graphics}
\usepackage{graphicx,psfrag}
\usepackage{graphics}
\usepackage{mathrsfs}
\usepackage{tabularx}
\usepackage{booktabs}
\usepackage{floatrow}
\usepackage{multirow}
\usepackage{graphicx}
\usepackage{caption}
\usepackage[ruled]{algorithm2e}
\newfloatcommand{capbtabbox}{table}[][\FBwidth]

\newtheorem{theo}{Theorem}[section]
\newtheorem{lem}[theo]{Lemma}

\newtheorem{cl}{Claim}[section]
\newtheorem{coro}[theo]{Corollary}

 \theoremstyle{definition}

\theoremstyle{remark}

\def\pf{\noindent{\bf Proof.\ }}
\def\qed{{\hfill\rule{4pt}{7pt}}}

\allowdisplaybreaks

\usepackage{indentfirst,subfig}
\begin{document}
    \captionsetup[figure]{labelfont={bf},name={Fig.},labelsep=period}

\begin{center} {\large Maximal and Maximum Dissociation Sets in General and Triangle-Free Graphs}
\end{center}
\pagestyle{empty}

\begin{center}
{
{\small Jianhua Tu$^{a,b,c}$, Yuxin Li$^{a}$, Junfeng Du$^{a,}$\footnote{Corresponding author.\\\indent \  E-mail addresses: tujh81@163.com (J. Tu), 2018200896@mail.buct.edu.cn (L. Zhang), dujf1990@163.com (J. Du)}}\\[2mm]
{\small $^a$ Department of mathematics, Beijing University of Chemical Technology, \\
\hspace*{1pt} Beijing, P.R. China 100029} \\[2mm]
{\small $^b$ Key Laboratory of Tibetan Information Processing and Machine Translation,\\
\hspace*{1pt} Qinghai Province, XiNing, P.R. China 810008}\\[2mm]
{\small $^c$ Key Laboratory of Tibetan Information Processing, Ministry of Education, \\
\hspace*{1pt} XiNing, P.R. China, 810008}\\
}
\end{center}
%
%
\begin{center}

\begin{abstract}
A subset of vertices $F$ in a graph $G$ is called a \emph{dissociation set} if the induced subgraph $G[F]$ of $G$ has maximum degree at most 1. A \emph{maximal dissociation set} of $G$ is a dissociation set which is not a proper subset of any other dissociation sets. A \emph{maximum dissociation set} is a dissociation set of maximum size. We show that every graph of order $n$ has at most $10^{\frac{n}{5}}$ maximal dissociation sets, and that every triangle-free graph of order $n$ has at most $6^{\frac{n}{4}}$ maximal dissociation sets. We also characterize the extremal graphs on which these upper bounds are attained. The tight upper bounds on the number of maximum dissociation sets in general and triangle-free graphs are also obtained.

\vskip 3mm


\end{abstract}
\end{center}

\baselineskip=0.24in
\section{Introduction}
Throughout this paper, we restrict our attention to simple, undirected and labeled graphs. Given a graph $G=(V,E)$ and a subset of vertices $S\subseteq V$, the subgraph of $G$ induced by $S$ is denoted by $G[S]$. We use $G-S$ to denote $G[V\setminus S]$, and $G-v$ for $G[V\setminus \{v\}]$. For $v\in V$, the neighborhood of $v$ in $G$ is denoted by $N_G(v)$, and the closed neighborhood is $N_G[v]=N_G(v)\cup \{v\}$. The degree of $v$ is denoted by $d_G(v)$, where $d_G(v)=|N_G(v)|$. A vertex of degree 1 is called a \emph{leaf}. We denote by $\delta(G)$ and $\Delta(G)$ the minimum and maximum degrees of the vertices of $G$.

Let $P_n$, $C_n$ and $K_n$ respectively denote the path, the cycle and the complete graph on $n$ vertices. $K_{m, n}$ denotes the complete bipartite graph with a partition $(X, Y)$, where $|X|=m$ and $|Y|=n$. Two graphs $G$ and $H$ are disjoint if $V(G)\cap V(H)=\emptyset$. For any two disjoint graphs $G$ and $H$, let $G\cup H$ denote the disjoint union of $G$ and $H$. For any positive integer $\ell$, $\ell G$ is the short notation for the disjoint union of $\ell$ copies of a graph $G$.

A subset of vertices $F$ in a graph $G$ is called a \emph{dissociation set} if the induced subgraph $G[F]$ has maximum degree at most 1. A \emph{maximal dissociation set} of $G$ is a dissociation set which is not a proper subset of any other dissociation sets. A maximum dissociation set is a dissociation set of maximum size. The concept of dissociation set was introduced by Yannakakis \cite{Yannakakis1981} and is a natural generalization of the well-known independent sets. The problem of finding a maximum dissociation set in a given graph was shown to be NP-complete even in planar line graph graphs of planar bipartite graphs \cite{Orlovich2011} and has been extensively studied \cite{Kardos2011, Orlovich2011, Xiao2017}.

In 1960s, Erd\H{o}s and Moser raised the problem of finding the maximum number of maximal cliques, or equivalently the maximum number of maximal independent sets among all graphs of order $n$. Since then the study of finding the maximum number of vertex subsets with given properties is an well-established area in graph theory and combinatorics. In particular, there has been a large number of results concerning the maximum number of maximal/maximum independent sets on many graph classes. For these results, we refer to \cite{Furedi1987,Griggs1988,Hujter1993,Jou2000,Koh2008,Liu1993,Mohr2018,Mohr2020,Moon1965,Sagan1988,Sagan2006,Wilf1986,Wloch2008,Zito1991}. Some other vertex subsets with given properties that have also been studied include minimal/minimum dominating sets \cite{Alvarado2019,Connolly2016,Fomin2005}, minimal connected vertex covers \cite{Golovach2018}, minimal feedback vertex sets \cite{Couturier2012,Fomin2008,Gaspers2013}, maximal induced matchings \cite{Basavaraju2014}, etc.

Inspired by these aforementioned results, Tu, Zhang and Shi \cite{Tu2019} considered the analogous problem for the \emph{maximum dissociation sets} and found the maximum number of maximum dissociation sets and the extremal graphs for trees with $n$ vertices. The set of all \emph{maximal dissociation sets} of a graph $G$ is denoted by $MD(G)$ and its cardinality by $\phi(G)$. In this paper, we study the number of \emph{maximal dissociation sets} in a graph of order $n$ and in a triangle-free graph of order $n$ and present the following two main theorems.


\begin{theo}\label{the1.1}
Let $G$ be a graph of order $n$. Then $\phi(G)\le 10^{\frac{n}{5}}$. Furthermore, $\phi(G)=10^{\frac{n}{5}}$ iff $n\equiv0\pmod5$ and $G$ is the disjoint union of $n/5$ graphs each of which is isomorphic to a graph obtained from $K_5$ by possibly deleting $0\leq i\leq 2$ non-adjacent edges.
\end{theo}

\begin{theo}\label{the1.2}
Let $G$ be a triangle-free graph of order $n$. Then $\phi(G)\le 6^{\frac{n}{4}}$. Furthermore, $\phi(G)= 6^{\frac{n}{4}}$ iff $n\equiv0\pmod4$ and $G\cong\frac{n}{4}C_4$.
\end{theo}

Obviously, the following corollary can be obtained by Theorem \ref{the1.2}.

\begin{coro}\label{cor1.1}
Let $G$ be a bipartite graph of order $n$. Then $\phi(G)\le 6^{\frac{n}{4}}$. Furthermore, $\phi(G)=6^{\frac{n}{4}}$ iff $n\equiv0\pmod4$ and $G\cong\frac{n}{4}C_4$.
\end{coro}

We let $\phi'(G)$ denote the number of maximum dissociation sets of $G$. Since every maximum dissociation set is also maximal, so $\phi'(G)\leq \phi(G)$. Given a family $\mathcal{G}$ of graphs on $n$ vertices, it is easy to see that if $G\in\mathcal{G}$ has the maximum number of maximal dissociation sets possible and every maximal dissociation set of $G$ is also a maximum dissociation set, then $G$ also has the maximum number of maximum dissociation sets possible. Since the graph $C_4$ and the graph obtained from $K_5$ by possibly deleting $0\le i\le 2$ non-adjacent edges satisfy the property that every maximal dissociation set is also a maximum dissociation set, we immediately have the following results.

\begin{theo}
Let $G$ be a graph of order $n$. Then $\phi'(G)\le 10^{\frac{n}{5}}$. Furthermore, $\phi'(G)=10^{\frac{n}{5}}$ iff $n\equiv0\pmod5$ and $G$ is the disjoint union of $n/5$ graphs each of which is isomorphic to a graph obtained from $K_5$ by possibly deleting $0\leq i\leq 2$ non-adjacent edges.
\end{theo}
\begin{theo}
Let $G$ be a triangle-free graph of order $n$. Then $\phi'(G)\le 6^{\frac{n}{4}}$. Furthermore, $\phi'(G)= 6^{\frac{n}{4}}$ iff $n\equiv0\pmod4$ and $G\cong\frac{n}{4}C_4$.
\end{theo}

The remainder of the paper is organized as follows. In Section 2, some preparatory lemmas are presented. In Section 3 we prove Theorem \ref{the1.1}. In Section 4 we prove Theorem \ref{the1.2}. In Section 5, we give the exact value of the maximum number of maximal dissociation sets in a graph of order $n$ and in a triangle-free graph of order $n$.

\section{Preparatory lemmas}
Given a graph $G$ and a vertex $v\in V(G)$, let
\begin{eqnarray*}
\phi(G,\overline{v}\ )&=&|\{F|F\in MD(G) \text{\ and }v\notin F\}|,\\
\phi(G,v^0)&=&|\{F|F\in MD(G), v\in F \text{\ and }d_{G[F]}(v)=0\}|,\\
\phi(G,v^1)&=&|\{F|F\in MD(G), v\in F \text{\ and }d_{G[F]}(v)=1\}|.
\end{eqnarray*}
For convenience, let $\alpha=10^{\frac{1}{5}}$ and $\beta=6^{\frac{1}{4}}$.

\begin{lem}\label{lem2.1}
For two disjoint graphs $G$ and $H$, \[\phi (G\cup H)=\phi (G)\cdot \phi (H).\]
\end{lem}

\begin{lem}\label{lem2.2}
Let $G$ be a graph and $v$ be a vertex of $G$,
$$\phi (G)\le \phi (G-v)+\phi (G-N[v])+\sum_{u\in N(v)}\phi (G-N[v]\cup N[u]).$$
Moreover, if there exists a vertex $w\in N(v)$ such that $N[w]\subseteq N[v]$, then
$$\phi (G)\le \phi (G-v)+\sum_{u\in N(v)}\phi (G-N[v]\cup N[u]).$$
\end{lem}

\pf Obviously,
\begin{eqnarray*}
\phi(G)&=&\phi(G,\overline{v})+\phi(G,v^0)+\phi(G,v^1)\\
&\leq&\phi (G-v)+\phi (G-N[v])+\sum_{u\in N(v)}\phi (G-N[v]\cup N[u]).
\end{eqnarray*}

If there exists a vertex $w\in N(v)$ such that $N[w]\subseteq N[v]$, then $\phi(G,v^0)=0$ and
\begin{eqnarray*}
\phi(G)&=&\phi(G,\overline{v})+\phi(G,v^1)\\
&\leq&\phi (G-v)+\sum_{u\in N(v)}\phi (G-N[v]\cup N[u]).
\end{eqnarray*}
 \qed

\begin{lem}\label{lem2.3}
Let $G$ be a graph. If $v$ is a leaf of $G$ and $N(v)=\{w\}$, then
$$\phi(G)\le \sum_{u\in N(w)\setminus\{v\}}\phi (G-N[w]\cup N[u])+\phi(G-v-w)+\phi(G-N[w]).$$
\end{lem}

\pf If a maximal dissociation set $F$ of $G$ doesn't contain the vertex $v$, then $w\in F$ and $d_{G[F]}(w)=1$. Thus,
\begin{eqnarray*}
\phi(G)&=&\phi(G,\overline{v})+\phi(G,v^0)+\phi(G,v^1)\\
&\le& \sum_{u\in N(w)\setminus\{v\}}\phi (G-N[w]\cup N[u])+\phi(G-v-w)+\phi(G-N[w]).
\end{eqnarray*}
\qed

\begin{lem}\label{lem2.4}
Let $G$ be a graph and $w$ be a vertex of $G$. If $w$ is adjacent to two leaves $v_1$ and $v_2$, then
$$\phi(G)\le \sum_{u\in N(w)\setminus\{v_1\}}\phi (G-N[w]\cup N[u])+\phi(G-\{w,v_1,v_2\})+\phi(G-N[w]).$$
\end{lem}

\pf Since $\phi(G-\{w,v_1\})=\phi(G-\{w,v_1,v_2\})$, we have
\begin{eqnarray*}
\phi(G)&=&\phi(G,\overline{v_1})+\phi(G,v_1^0)+\phi(G,v_1^1)\\
&\leq&\sum_{u\in N(w)\setminus\{v_1\}}\phi (G-N[w]\cup N[u])+\phi(G-\{w,v_1,v_2\})+\phi(G-N[w]).
\end{eqnarray*}\qed

\begin{lem}\label{lem2.5}
For a path $P_n$, \[\phi(P_n)<0.81\beta^n.\]
\end{lem}

\pf We prove the lemma by induction on $n$. When $n\le 4$, it is easy to check that the result holds. We assume that $n\geq5$ and the result holds for any path with at most $n-1$ vertices. Let $v$ be a leaf of $P_n$, $w$ be the neighbour of $v$, $u$ be another neighbour of $w$. By Lemma \ref{lem2.3} and the inductive hypothesis, we have
\begin{eqnarray*}
\phi(P_n) & \le & \phi(P_n-N[w]\cup N[u])+\phi(P_n-v-w)+\phi(P_n-N[w]) \\
& < & 0.81\beta^{n-4}+0.81\beta^{n-2}+0.81\beta^{n-3} \\
& = & 0.81\beta^n(\beta^{-4}+\beta^{-2}+\beta^{-3})<0.81\beta^n.
\end{eqnarray*} \qed

\begin{lem}\label{lem2.6}
For a cycle $C_n$, \[\phi(C_n)\le\beta^n.\]
Furthermore, $\phi(C_n)=\beta^n$ iff $n=4$.
\end{lem}

\pf It is easy to check that the result holds for $n=3$ or $n=4$. Thus, we assume that $n\geq 5$.
Let $v$ be a vertex of $C_n$. By Lemma \ref{lem2.2} and Lemma \ref{lem2.5}, we have
\begin{eqnarray*}
\phi(C_n) & \le & \phi(C_n-v)+\phi(C_n-N[v])+\sum_{u\in N(v)}\phi(C_n-N[v]\cup N[u]) \\
& < & 0.81\beta^{n-1}+0.81\beta^{n-3}+2\cdot 0.81\beta^{n-4} \\
& = & 0.81\beta^n(\beta^{-1}+\beta^{-3}+2\beta^{-4})<\beta^n.
\end{eqnarray*} \qed

\section{Proof of Theorem \ref{the1.1}}

We prove the theorem by induction on $n$. It is easy to check that the theorem is true for $n\le5$. Thus, we assume that $n\ge6$ and the theorem is true for any graph with at most $n-1$ vertices.

Suppose that $G$ is disconnected. Let $G_1$ be a component of $G$, $G_2$ be the union of other components of $G$, $n_1=|V(G_1)|$ and $n_2=|V(G_2)|$. By the inductive hypothesis,
$$\phi(G)= \phi(G_1)\cdot \phi(G_2)\le \alpha^{n_1}\cdot \alpha^{n_2}=\alpha^n,$$
moreover, $\phi(G)=\alpha^n$ implies that $n_1\equiv 0\pmod 5$, $n_2\equiv 0\pmod5$, both $G_1$ and $G_2$ are the disjoint union of some graphs each of which is isomorphic to a graph obtained from $K_5$ by possibly deleting $0\le i\le 2$ non-adjacent edges. It follows that $\phi(G)=\alpha^n$ iff $n\equiv 0\pmod5$ and $G$ is the disjoint union of $n/5$ graphs each of which is isomorphic to a graph obtained from $K_5$ by possibly deleting $0\le i\le 2$ non-adjacent edges.

Now, suppose that $G$ is connected. We first present the following two claims.

\begin{cl}\label{cla3.1}
Let $v$ be a leaf of $G$ and $w$ be the neighbor of the vertex $v$. If $d(w)=2$ and another neighbor of $w$ has at least two neighbors, then \[\phi(G)\le \alpha^{n-4}+\alpha^{n-2}+\alpha^{n-3}.\]
\end{cl}
\pf Let $u$ be another neighbor of the vertex $w$. By Lemma \ref{lem2.3} and the inductive hypothesis, we have
\begin{eqnarray*}
\phi(G) & \le & \phi (G-N[w]\cup N[u])+\phi(G-v-w)+\phi(G-N[w]) \\
& \le & \alpha^{n-4}+\alpha^{n-2}+\alpha^{n-3}.\\
\end{eqnarray*}\qed

\begin{figure}[htb]
     \centering
     \includegraphics[scale=0.5]{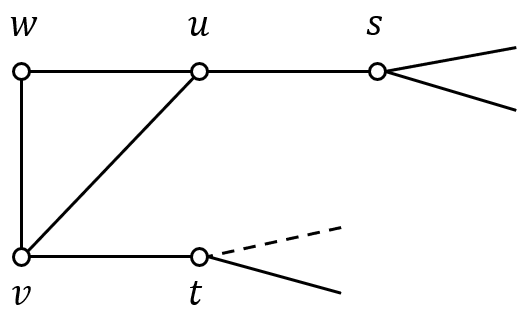}
     \caption{}
     \label{fig1}
 \end{figure}

\begin{cl}\label{cla3.2}
If $G$ contains a substructure shown in Fig.\ref{fig1}, where $d(v)=d(u)=d(s)=3$, $d(w)=2$ and $d(t)\ge2$, then $\phi(G)<0.832\alpha^n$.
\end{cl}

 \pf Similarly, by Lemma \ref{lem2.3} and the inductive hypothesis, we have \[\phi(G-v)\le \alpha^{n-6}+\alpha^{n-3}+\alpha^{n-4}.\]
By Lemma \ref{lem2.2} and the inductive hypothesis,
\begin{eqnarray*}
\phi(G) & \le & \phi(G-v)+\sum_{t\in N(v)}\phi(G-N[v]\cup N[t]) \\
& \le & \alpha^{n-6}+\alpha^{n-3}+\alpha^{n-4}+\alpha^{n-4}+2\alpha^{n-5} \\
& = & \alpha^n(\alpha^{-6}+\alpha^{-3}+2\alpha^{-4}+2\alpha^{-5})<0.832\alpha^n.
\end{eqnarray*}
\qed

It is worth mentioning that the two claims hold regardless of whether $G$ is connected or not. We proceed to prove that $\phi(G)<\alpha^n$ and distinguish the following three cases.

\textbf{Case 1. } $\Delta(G)=2$.

In this case, $G$ is either a path or a cycle. By Lemma \ref{lem2.5} and \ref{lem2.6}, we have $\phi(G)<\alpha^n$.

\textbf{Case 2. } $\Delta(G)=3$.

\textit{Subcase 2.1. } $\delta(G)=1$.

Let $v$ be a leaf of $G$, $N(v)=\{w\}$. If $d(w)=2$, then by Claim \ref{cla3.1}, $\phi(G)\le \alpha^{n-4}+\alpha^{n-2}+\alpha^{n-3}<\alpha^n$. If $d(w)=3$, then by Lemma \ref{lem2.3} and the inductive hypothesis,
\begin{eqnarray*}
\phi(G) & \le & \sum_{u\in N(w)\setminus\{v\}}\phi (G-N[w]\cup N[u])+\phi(G-v-w)+\phi(G-N[w]) \\
& \le & \alpha^{n-5}+\alpha^{n-4}+\alpha^{n-2}+\alpha^{n-4} \\
& = & \alpha^n(\alpha^{-5}+2\alpha^{-4}+\alpha^{-2}) <  \alpha^n.
\end{eqnarray*}

\textit{Subcase 2.2. } $\delta(G)=2$.


\textit{Subcase 2.2.1.} There exists a triangle $G[\{w,u,v\}]\cong K_3$ in $G$ such that $d(w)=2$ and $d(v)=3$.

If $d(u)=2$, then by Lemma \ref{lem2.1}, \ref{lem2.2} and the inductive hypothesis, we have
\begin{eqnarray*}
\phi(G) & \le & \phi(G-v)+\sum_{t\in N(v)}\phi(G-N[v]\cup N[t]) \\
& = & \phi(G-\{v,u,w\})+\sum_{t\in N(v)}\phi(G-N[v]\cup N[t]) \\
& \le & \alpha^{n-3}+2\alpha^{n-4}+\alpha^{n-5}\\
&<& \alpha^n(\alpha^{-3}+2\alpha^{-4}+\alpha^{-5})< \alpha^n.
\end{eqnarray*}
If $d(u)=3$, then by Claim \ref{cla3.1}, $\phi(G-v)\le  \alpha^{n-5}+\alpha^{n-3}+\alpha^{n-4}.$
Thus, by Lemma \ref{lem2.2} and the inductive hypothesis,
\begin{eqnarray*}
\phi(G) & \le & \phi(G-v)+\sum_{t\in N(v)}\phi(G-N[v]\cup N[t]) \\
& \le & \alpha^{n-5}+\alpha^{n-3}+\alpha^{n-4}+\alpha^{n-5}+2\alpha^{n-4} \\
& = & \alpha^n(2\alpha^{-5}+\alpha^{-3}+3\alpha^{-4})<\alpha^n.
\end{eqnarray*}

\textit{Subcase 2.2.2}. For any triangle of $G$, all its vertices have degree 3.

Let $w$ be a vertex of degree 2 of $G$ that is adjacent to a vertex $v$ of degree 3. Let $u$ be another neighbor of the vertex $w$. Clearly, $uv\notin E(G)$.

When $d(u)=3$, $\phi(G-v)\le  2\alpha^{n-5}+\alpha^{n-6}+\alpha^{n-3}$ and
\begin{eqnarray*}
\phi(G) & \le & \phi(G-v)+\phi(G-N[v])+\sum_{t\in N(v)}\phi(G-N[v]\cup N[t]) \\
& \le & 2\alpha^{n-5}+\alpha^{n-6}+\alpha^{n-3}+\alpha^{n-4}+3\alpha^{n-5} \\
& = & \alpha^n(5\alpha^{-5}+\alpha^{-3}+\alpha^{-4}+\alpha^{-6})<\alpha^n.
\end{eqnarray*}

Assume that $d(u)=2$. Let $s$ be another neighbor of $u$. If $sv\in E(G)$ and $d(s)=2$, then $G-v=P_3\cup (G-\{v,w,u,s\})$, which implies that $\phi(G-v)\le 3\alpha^{n-4}$. By Lemma \ref{lem2.2} and the inductive hypothesis, we have
\begin{eqnarray*}
\phi(G) & \le & \phi(G-v)+\phi(G-N[v])+\sum_{t\in N(v)}\phi(G-N[v]\cup N[t]) \\
& \le & 3\alpha^{n-4}+\alpha^{n-4}+3\alpha^{n-5} \\
& = & \alpha^n(4\alpha^{-4}+3\alpha^{-5})<\alpha^n.
\end{eqnarray*}
Otherwise, either $sv\notin E(G)$ or $d(s)=3$. It is easy to show that $\phi(G-v)\le \alpha^{n-5}+\alpha^{n-3}+\alpha^{n-4}$. Furthermore, since for any triangle of $G$, all its vertices have degree 3, we have
\begin{eqnarray*}
\phi(G) & \le & \phi(G-v)+\phi(G-N[v])+\sum_{t\in N(v)}\phi(G-N[v]\cup N[t]) \\
& \le & \alpha^{n-5}+\alpha^{n-3}+\alpha^{n-4}+\alpha^{n-4}+3\alpha^{n-5} \\
& = & \alpha^n(4\alpha^{-5}+\alpha^{-3}+2\alpha^{-4})<\alpha^n.
\end{eqnarray*}

\textit{Subcase 2.3. } $\delta(G)=3$.

In this subcase, $G$ is a cubic graph. Let $v$ be a vertex in $G$. There are three possibilities for $G[N(v)]$:
either $G[N(v)]\cong \overline{K_3}$, or $G[N(v)]\cong P_3$, or $G[N(v)]\cong P_2\cup K_1$.
If $G[N(v)]\cong \overline{K_3}$, then by Lemma \ref{lem2.2} and the inductive hypothesis, we have
\begin{eqnarray*}
\phi(G) & \le & \phi(G-v)+\phi(G-N[v])+\sum_{u\in N(v)}\phi(G-N[v]\cup N[u]) \\
& \le & \alpha^{n-1}+\alpha^{n-4}+3\alpha^{n-6}\\
& = & \alpha^n(\alpha^{-1}+\alpha^{-4}+3\alpha^{-6})<\alpha^n.
\end{eqnarray*}

If $G[N(v)]\cong P_3$, then by Lemma \ref{lem2.2} and the inductive hypothesis, we have
\begin{eqnarray*}
\phi(G) & \le & \phi(G-v)+\sum_{u\in N(v)}\phi(G-N[v]\cup N[u]) \\
& \le & \alpha^{n-1}+\alpha^{n-4}+2\alpha^{n-5} \\
& = & \alpha^n(\alpha^{-1}+\alpha^{-4}+2\alpha^{-5})<\alpha^n.
\end{eqnarray*}

\begin{figure}[htb]
     \centering
     \includegraphics[scale=0.5]{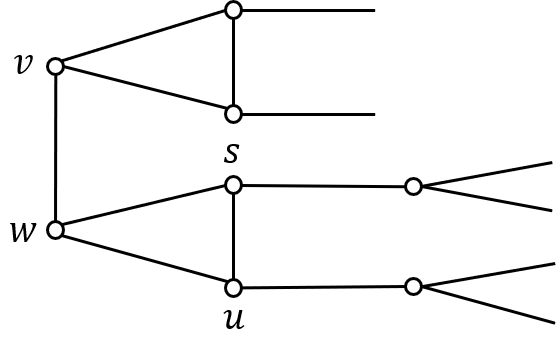}
     \caption{}
     \label{fig2}
 \end{figure}

Now, consider the last possibility, i.e., $G[N(v)]\cong P_2\cup K_1$.  Let $w$ be the isolated vertex in $G[N(v)]$. By the aforementioned possibilities, we can assume that for any vertex $u\in V(G)$, $G[N(u)]\cong P_2\cup K_1$, otherwise one can deduce that $\phi(G)<\alpha^n$. If $n=6$, then it is easy to check that $\phi(G)=9<\alpha^6$. If $n\ge7$, $G$ contains a substructure shown in Fig.\ref{fig2}, furthermore, $G-v$ contains a substructure shown in Fig.\ref{fig1}.
By Claim \ref{cla3.2}, $\phi(G-v)<0.832\alpha^{n-1}$. By Lemma \ref{lem2.2} and the inductive hypothesis,
\begin{eqnarray*}
\phi(G) & \le & \phi(G-v)+\phi(G-N[v])+\sum_{t\in N(v)}\phi(G-N[v]\cup N[t]) \\
& \le & 0.832\alpha^{n-1}+\alpha^{n-4}+2\alpha^{n-5}+\alpha^{n-6} \\
& = & \alpha^n(0.832\alpha^{-1}+\alpha^{-4}+2\alpha^{-5}+\alpha^{-6})<\alpha^n.
\end{eqnarray*}

\textbf{Case 3. } $\Delta(G) \ge 4$.

Let $v$ be a vertex of degree $\Delta$. If for any $u\in N(v)$, $N[u]\not\subseteq N[v]$, then by Lemma \ref{lem2.2} and the inductive hypothesis, we have
\begin{eqnarray*}
\phi(G) & \le & \phi(G-v)+\phi(G-N[v])+\sum_{u\in N(v)}\phi(G-N[v]\cup N[u]) \\
& \le & \alpha^{n-1}+\alpha^{n-\Delta-1}+\Delta \alpha^{n-\Delta-2} \\
& = & \alpha^n(\alpha^{-1}+\alpha^{-\Delta-1}+\Delta \alpha^{-\Delta-2}) <\alpha^n.
\end{eqnarray*}
Otherwise, there exists a vertex $t\in N(v)$ such that $N[t]\subseteq N[v]$. By Lemma \ref{lem2.2} and the inductive hypothesis,
\begin{eqnarray*}
\phi(G) & \le & \phi(G-v)+\sum_{u\in N(v)}\phi(G-N[v]\cup N[u]) \\
& \le & \alpha^{n-1}+\Delta \alpha^{n-\Delta-1} \\
& = & \alpha^n(\alpha^{-1}+\Delta \alpha^{-\Delta-1})
\end{eqnarray*}
When $\Delta>4$, $\alpha^{-1}+\Delta \alpha^{-\Delta-1}<1$, it follow that $\phi(G)<\alpha^n$. When $\Delta=4$, since $n\ge6$, there exists a vertex $w\in N(v)$ such that $N[w]\not\subseteq N[v]$. At this moment, we have
\begin{eqnarray*}
\phi(G) & \le & \phi(G-v)+\sum_{u\in N(v)}\phi(G-N[v]\cup N[u]) \\
& \le & \alpha^{n-1}+3\alpha^{n-5}+\alpha^{n-6} \\
& = & \alpha^n(\alpha^{-1}+3\alpha^{-5}+\alpha^{-6})<\alpha^n.
\end{eqnarray*}

We complete the proof of the theorem.
\qed

\section{Proof of Theorem \ref{the1.2}}

We prove the theorem by induction on $n$.  It is easy to check that the theorem is true for $n\le4$. Thus, we assume that $n\ge5$ and the theorem is true for any graph with at most $n-1$ vertices.

Suppose that $G$ is disconnected. Let $G_1$ be a component of $G$, $G_2$ be the union of other components of $G$, $n_1=|V(G_1)|$ and $n_2=|V(G_2)|$. By the inductive hypothesis,
$$\phi(G)= \phi(G_1)\cdot \phi(G_2)\le \beta^{n_1}\cdot \beta^{n_2}=\beta^n,$$
moreover, $\phi(G)=\beta^n$ implies that $n_1\equiv 0\pmod 4$, $n_2\equiv 0\pmod 4$, both $G_1$ and $G_2$ are the disjoint union of some graphs each of which is isomorphic to $C_4$. It follows that $\phi(G)=\beta^n$ iff $n\equiv 0\pmod4$ and $G\cong\frac{n}{4}C_4$.

Now, suppose that $G$ is connected.

\begin{cl}\label{cl4-1}
If $G$ contains a leaf, then $\phi(G)<0.93\beta^n$.
\end{cl}

\pf Let $v$ be a leaf of $G$. Let $N(v)=\{w\}$ and $d(w)=d$. By Lemma \ref{lem2.3}, we have
\begin{eqnarray*}
\phi(G) & \le & \sum_{u\in N(w)\setminus \{v\}}\phi (G-N[w]\cup N[u])+\phi(G-v-w)+\phi(G-N[w]) \\
& \le & (d-1)\cdot \beta^{n-d-1}+\beta^{n-2}+\beta^{n-d-1} \\
& = & \beta^n(\beta^{-2}+d\cdot \beta^{-d-1}).
\end{eqnarray*}
 When $d\ge 1$, $\phi(G)<0.93\beta^n$. \qed

It is worth mentioning that the claim holds regardless of whether $G$ is connected
or not. By Claim \ref{cl4-1}, we assume that $\delta(G)\ge2$. We proceed to prove that $\phi(G)<\beta^n$ and distinguish the following four cases.

\textbf{Case 1.} $\Delta(G)=2$.

In this case, $G$ is a cycle. By Lemma \ref{lem2.6}, we have $\phi(G)<\beta^n$.

\textbf{Case 2.} $\Delta(G)=3$.

\textit{Subcase 2.1.} $\delta(G)=2$.

In this subcase, there must be a vertex of degree 2, say $w$, that is adjacent to a vertex $v$ of degree 3. Let $N(w)=\{u,v\}$.

When $d(u)=2$, let $N(u)=\{w, s\}$. If $sv\in E(G)$ and $d(s)=2$, then $G-v\cong P_3\cup (G-\{v, w, u, s\})$, which implies that $\phi(G-v)\le 3\beta^{n-4}$ and $\phi(G-N[v])  \le \beta^{n-5}$. By Lemma \ref{lem2.2} and the inductive hypothesis,
\begin{eqnarray*}
\phi(G) & \le & \phi(G-v)+\phi(G-N[v])+\sum_{t\in N(v)}\phi(G-N[v]\cup N[t]) \\
& \le & 3\beta^{n-4}+\beta^{n-5}+3\beta^{n-5} \\
& = & \beta^n(3\beta^{-4}+4\beta^{-5})<\beta^n.
\end{eqnarray*}
If $sv\in E(G)$ and $d(s)=3$, then $G-v$ contains a leaf, $\phi(G-v) \le \beta^{n-5}+\beta^{n-3}+\beta^{n-4}$ and $\phi(G-N[v])  \le \beta^{n-5}$.
%
By Lemma \ref{lem2.2} and the inductive hypothesis,
\begin{eqnarray*}
\phi(G) & \le & \phi(G-v)+\phi(G-N[v])+\sum_{t\in N(v)}\phi(G-N[v]\cup N[t]) \\
& \le & \beta^{n-5}+\beta^{n-3}+\beta^{n-4}+\beta^{n-5}+2\beta^{n-5}+\beta^{n-6} \\
& = & \beta^n(4\beta^{-5}+\beta^{-3}+\beta^{-4}+\beta^{-6})<\beta^n.
\end{eqnarray*}
Finally, if $sv\notin E(G)$, then $G-v$ contains a leaf and
$\phi(G-v)\le \beta^{n-5}+\beta^{n-3}+\beta^{n-4}$. Furthermore, through a careful analysis, we have
\begin{eqnarray*}
\phi(G-N[v])&\le&\max\{\beta^{n-6},3\beta^{n-7},\beta^{n-8}+\beta^{n-6}+\beta^{n-7},4\beta^{n-8},\beta^{n-9}+2\beta^{n-8}+\beta^{n-6}\}\\
&\le&\beta^{n-9}+2\beta^{n-8}+\beta^{n-6}.
\end{eqnarray*}
Since $G$ is a triangle-free graph,
\begin{eqnarray*}
\phi(G) & \le & \phi(G-v)+\phi(G-N[v])+\sum_{t\in N(v)}\phi(G-N[v]\cup N[t]) \\
& \le &\beta^{n-5}+\beta^{n-3}+\beta^{n-4}+\beta^{n-9}+2\beta^{n-8}+\beta^{n-6}+3\beta^{n-5} \\
& = & \beta^n(4\beta^{-5}+\beta^{-3}+\beta^{-4}+\beta^{-9}+2\beta^{-8}+\beta^{-6})<\beta^n.
\end{eqnarray*}

When $d(u)=3$, if $n=5$, then $G\cong K_{2, 3}$ and $\phi(K_{2, 3})=8<\beta^5$. If $n\ge 6$, then $G-v$ contains at least a leaf and $\phi(G-v)\le 2\beta^{n-6}+\beta^{n-3}+\beta^{n-5}$. Thus, by Lemma \ref{lem2.2} and the inductive hypothesis,
\begin{eqnarray*}
\phi(G) & \le & \phi(G-v)+\phi(G-N[v])+\sum_{t\in N(v)}\phi(G-N[v]\cup N[t]) \\
& \le & 2\beta^{n-6}+\beta^{n-3}+\beta^{n-5}+\beta^{n-4}+3\beta^{n-5} \\
& = & \beta^n(2\beta^{-6}+\beta^{-3}+\beta^{-4}+4\beta^{-5})<\beta^n.
\end{eqnarray*}

\textit{subcase 2.2.} $\delta(G)=3$.

\begin{figure}[htb]
     \centering
     \includegraphics[scale=0.5]{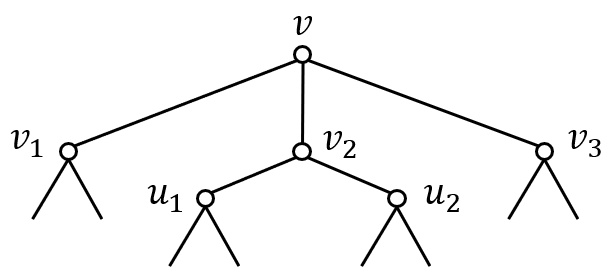}
     \caption{}
     \label{fig4}
 \end{figure}

In this subcase, $G$ is a cubic graph and contains a substructure shown in Fig.\ref{fig4}.
%
By Lemma \ref{lem2.2} and the inductive hypothesis,
\begin{eqnarray*}
\phi(G) & \le & \phi(G-v)+\phi(G-N[v])+\sum_{u\in N(v)}\phi(G-N[v]\cup N[u]) \\
& \le & \phi(G-v)+\beta^{n-4}+3\beta^{n-6}.
\end{eqnarray*}

\textit{Subcase 2.2.1. } $|N(u_1)\cap N(u_2)|=3$.

\begin{figure}[htb]
     \centering
     \includegraphics[scale=0.5]{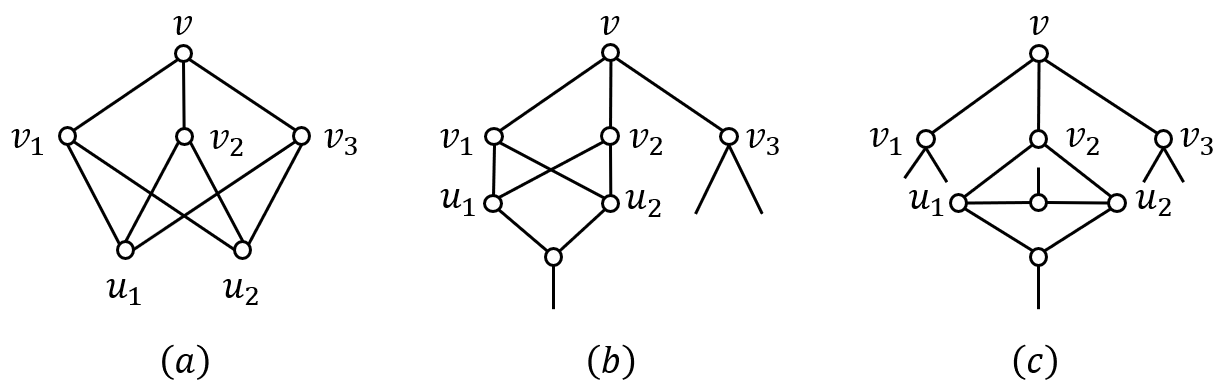}
     \caption{Three substructures}
     \label{fig5}
 \end{figure}

Now, $G$ contains one of three substructures shown in Fig.\ref{fig5}.

If $G$ contains a substructure shown in Fig.\ref{fig5} (a), then $G\cong K_{3,3}$ and $\phi(G)=11<\beta^6$.

If $G$ contains a substructure shown in Fig.\ref{fig5} (b), then
\begin{eqnarray*}
\phi(G-v) & \le & \phi(G-v-u_1)+\beta^{n-5}+2\beta^{n-6}+\beta^{n-7},\\
\phi(G-v-u_1) & \le & 2\beta^{n-6}+\beta^{n-7}+\beta^{n-5}.
\end{eqnarray*}
Hence, $\phi(G)  \le  \beta^n(2\beta^{-5}+2\beta^{-7}+7\beta^{-6}+\beta^{-4})<\beta^n$.

If $G$ contains a substructure shown in Fig.\ref{fig5} (c), then
\begin{eqnarray*}
\phi(G-v) & \le & \phi(G-v-u_1)+\beta^{n-5}+\beta^{n-6}+2\beta^{n-7},\\
\phi(G-v-u_1) & \le & 2\beta^{n-7}+\beta^{n-4}+\beta^{n-6}.
\end{eqnarray*}
Hence, $\phi(G)\le\beta^n(\beta^{-5}+4\beta^{-7}+5\beta^{-6}+2\beta^{-4})<\beta^n.$

\textit{Subcase 2.2.2 } $|N(u_1)\cap N(u_2)|=2.$

\begin{figure}[htb]
     \centering
     \includegraphics[scale=0.5]{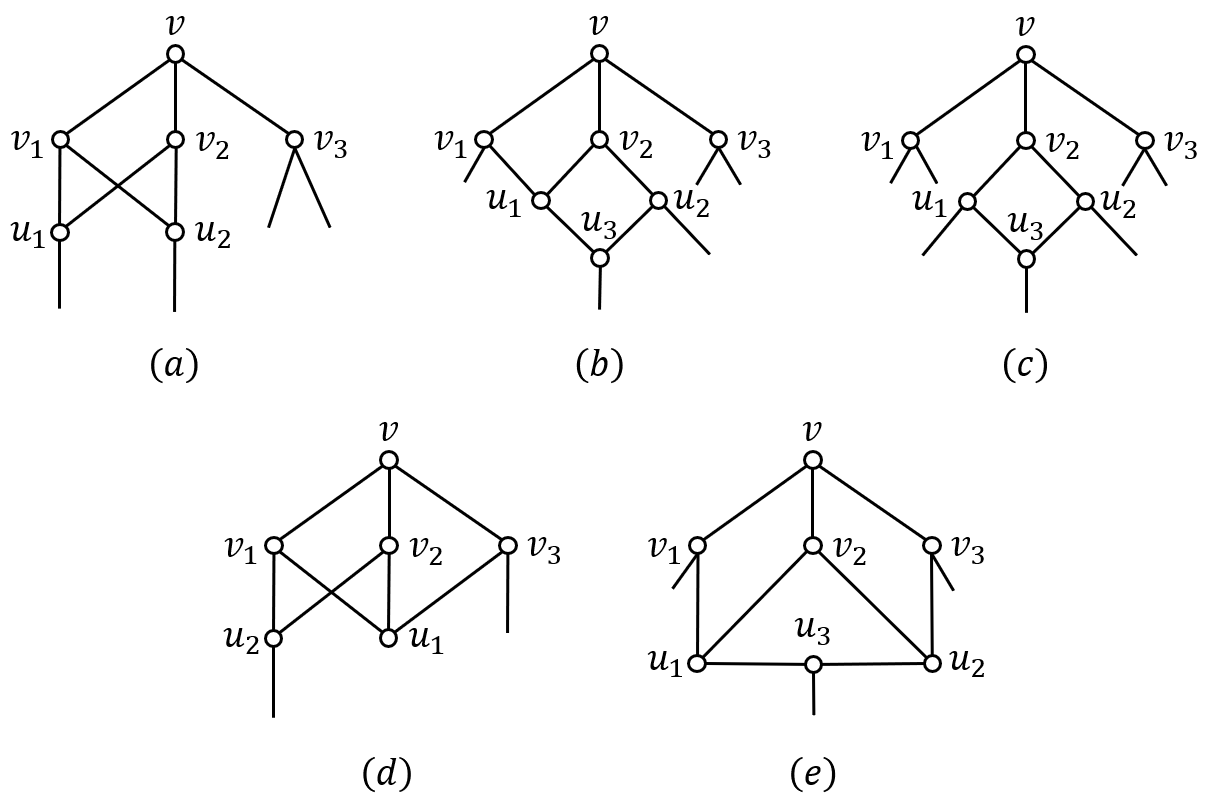}
     \caption{Five substructures}
     \label{fig6}
 \end{figure}

Now, $G$ contains one of five substructures shown in Fig.\ref{fig6}
.

If $G$ contains a substructure shown in Fig.\ref{fig6} $(a)$, or in Fig.\ref{fig6} $(b)$ or in Fig.\ref{fig6} $(c)$, then
\begin{eqnarray*}
\phi(G-v) & \le & \phi(G-v-u_1)+\beta^{n-5}+2\beta^{n-6}+\beta^{n-7},\\
\phi(G-v-u_1) & \le & \beta^{n-7}+\beta^{n-8}+\beta^{n-4}+\beta^{n-6}.
\end{eqnarray*}
Hence, $\phi(G) \le \beta^n (\beta^{-5}+2\beta^{-7}+6\beta^{-6}+2\beta^{-4}+\beta^{-8})<\beta^n$.

If $G$ contains a substructure shown in Fig.\ref{fig6} $(d)$, then
\begin{eqnarray*}
\phi(G-v) & \le & \phi(G-v-u_1)+\beta^{n-5}+3\beta^{n-6},\\
\phi(G-v-u_1) & \le & 2\beta^{n-6}+\beta^{n-8}+\beta^{n-5}.
\end{eqnarray*}
Hence, $\phi(G) \le \beta^n(2\beta^{-5}+8\beta^{-6}+\beta^{-4}+\beta^{-8})<\beta^n$.

If $G$ contains a substructure shown in Fig.\ref{fig6} $(e)$, then
\begin{eqnarray*}
\phi(G-v) & \le & \phi(G-v-u_1)+\beta^{n-5}+2\beta^{n-6}+\beta^{n-7},\\
\phi(G-v-u_1) & \le & 2\beta^{n-7}+\beta^{n-4}+\beta^{n-6}.
\end{eqnarray*}
Hence, $\phi(G) \le \beta^n(\beta^{-5}+3\beta^{-7}+6\beta^{-6}+2\beta^{-4})<\beta^n$.

\textit{Subcase 2.2.3. } $|N(u_1)\cap N(u_2)|=1$.

\begin{figure}[htb]
     \centering
     \includegraphics[scale=0.5]{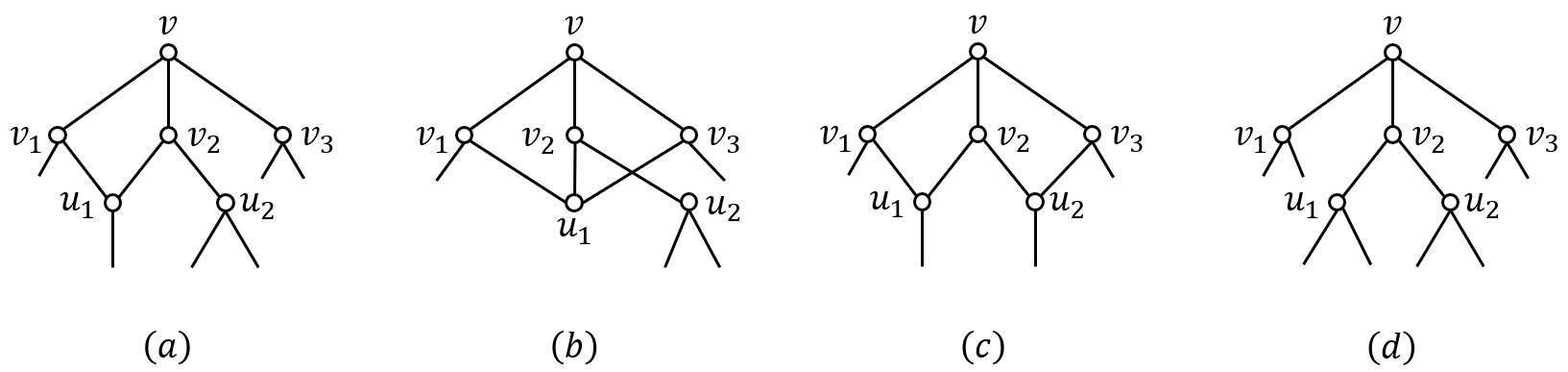}
     \caption{Four substructures}
     \label{fig8}
 \end{figure}

Now, $G$ contains one of four substructures shown in Fig.\ref{fig8}.

If $G$ contains a substructure shown in Fig.\ref{fig8} (a), then
\begin{eqnarray*}
\phi(G-v) & \le & \phi(G-v-u_1)+\beta^{n-5}+2\beta^{n-6}+\beta^{n-7},\\
\phi(G-v-u_1) & \le & 2\beta^{n-8}+\beta^{n-4}+\beta^{n-6}.
\end{eqnarray*}
Hence, $\phi(G)\le\beta^n (\beta^{-5}+\beta^{-7}+6\beta^{-6}+2\beta^{-4}+2\beta^{-8})<\beta^n.$

If $G$ contains a substructure shown in Fig.\ref{fig8} (b), then
\begin{eqnarray*}
\phi(G-v) & \le & \phi(G-v-u_1)+\beta^{n-5}+3\beta^{n-6},\\
\phi(G-v-u_1) & \le & 2\beta^{n-8}+\beta^{n-4}+\beta^{n-6}.
\end{eqnarray*}
Hence, $\phi(G)\le \beta^n(\beta^{-5}+7\beta^{-6}+2\beta^{-4}+2\beta^{-8})<\beta^n$.

If $G$ contains a substructure shown in Fig.\ref{fig8} (c), then
\begin{eqnarray*}
\phi(G-v) & \le & \phi(G-v-u_1)+\beta^{n-5}+2\beta^{n-6}+\beta^{n-7},\\
\phi(G-v-u_1) & \le & \beta^{n-7}+\beta^{n-8}+\beta^{n-4}+\beta^{n-6}.
\end{eqnarray*}
Hence, $\phi(G)\le \beta^n(\beta^{-5}+2\beta^{-7}+6\beta^{-6}+2\beta^{-4}+\beta^{-8})<\beta^n$.

If $G$ contains a substructure shown in Fig.\ref{fig8} (d), then
\begin{eqnarray*}
\phi(G-v) & \le & \phi(G-v-u_1)+\beta^{n-5}+\beta^{n-6}+2\beta^{n-7},\\
\phi(G-v-u_1) & \le & 2\beta^{n-8}+\beta^{n-4}+\beta^{n-6}.
\end{eqnarray*}
Hence, $\phi(G)\le \beta^n(\beta^{-5}+2\beta^{-7}+5\beta^{-6}+2\beta^{-4}+2\beta^{-8})<\beta^n$.

\textbf{Case 3. } $\Delta(G)=4$.

Let $v$ be a vertex of degree 4. If there exists a vertex $w$ of degree 2 in $N(v)$, then $G-v$ contains a leaf and by Claim \ref{cl4-1}, $\phi(G-v)<0.93\beta^{n-1}$. By Lemma \ref{lem2.2} and the inductive hypothesis,
\begin{eqnarray*}
\phi(G) & \le & \phi(G-v)+\phi(G-N[v])+\sum_{u\in N(v)}\phi(G-N[v]\cup N[u]) \\
& \le & 0.93\beta^{n-1}+\beta^{n-5}+4\beta^{n-6}\\
& = & \beta^n(0.93\beta^{-1}+\beta^{-5}+4\beta^{-6})<\beta^n.
\end{eqnarray*}
If for any $u\in N(v)$, $d(u)>2$. By Lemma \ref{lem2.2} and the inductive hypothesis,
\begin{eqnarray*}
\phi(G) & \le & \phi(G-v)+\phi(G-N[v])+\sum_{u\in N(v)}\phi(G-N[v]\cup N[u]) \\
& \le & \beta^{n-1}+\beta^{n-5}+4\beta^{n-7}\\
& = & \beta^n(\beta^{-1}+\beta^{-5}+4\beta^{-7})<\beta^n.
\end{eqnarray*}

\textbf{Case 4. } $\Delta(G) \ge 5$.

Let $v$ be a vertex of degree $\Delta$. By Lemma \ref{lem2.2} and the inductive hypothesis,
\begin{eqnarray*}
\phi(G) & \le & \phi(G-v)+\phi(G-N[v])+\sum_{u\in N(v)}\phi(G-N[v]\cup N[u]) \\
& \le & \beta^{n-1}+\beta^{n-\Delta-1}+\Delta \beta^{n-\Delta-2}\\
& = & \beta^n(\beta^{-1}+\beta^{-\Delta-1}+\Delta \beta^{n-\Delta-2}).
\end{eqnarray*}
When $\Delta \ge 5$, $\beta^{-1}+\beta^{-\Delta-1}+\Delta \beta^{n-\Delta-2}<1$, so $\phi(G)<\beta^n$.

We complete the proof of the theorem.\qed


%

\section{Concluding remarks}

In this paper, we present the upper bounds on the number of maximal dissociation sets in a general graph of order $n$ and in a triangle-free graphs of order $n$. In fact, we can give the exact value of the maximum number of maximal dissociation sets in a general graph of order $n$ and in a triangle-free graphs of order $n$, however, the proofs are similar and tedious, so we omit them.

\begin{theo}
If $G$ is a graph with $n\ge 8$ vertices, then
\begin{equation*}
\phi(G)\le \left\{
\begin{aligned}
&10^t, & n=5t \quad \ \  \\
&15\cdot 10^{t-1}, & n=5t+1 \\
&225\cdot 10^{t-2}, & n=5t+2 \\
&36\cdot 10^{t-1}, & n=5t+3 \\
&6\cdot 10^t, & n=5t+4
\end{aligned}
\right.
\end{equation*}
where $t$ is a positive integer, and the equality holds iff
\begin{equation*}
G \cong \left\{
\begin{aligned}
&tK_5^*, & n=5t \quad \ \  \\
&K_6^*\cup (t-1)K_5^*, & n=5t+1 \\
&2K_6^*\cup (t-2)K_5^*, & n=5t+2 \\
&2K_4^*\cup (t-1)K_5^*, & n=5t+3 \\
&K_4^*\cup tK_5^*, & n=5t+4
\end{aligned}
\right.
\end{equation*}
where $K_m^*$ is obtained from $K_m$ by possibly deleting $0\leq i\leq m/2$ non-adjacent edges.
\end{theo}

\begin{theo}
If $G$ is a triangle-free graph with $n\ge 4$ vertices, then
\begin{equation*}
\phi(G)\le \left\{
\begin{aligned}
&6^t, & n=4t \quad \ \  \\
&8\cdot 6^{t-1}, & n=4t+1 \\
&11\cdot 6^{t-1}, & n=4t+2 \\
&3\cdot 6^t, & n=4t+3
\end{aligned}
\right.
\end{equation*}
where $t$ is a positive integer, and the equality holds iff
\begin{equation*}
G \cong \left\{
\begin{aligned}
&tC_4, & n=4t \quad \ \  \\
&K_{2, 3}\cup (t-1)C_4, & n=4t+1 \\
&K_{3, 3}\cup (t-1)C_4, & n=4t+2 \\
&P_3\cup tC_4, & n=4t+3
\end{aligned}
\right.
\end{equation*}
\end{theo}

\bibliographystyle{unsrt}

\end{document}